\documentclass[11pt]{article}
\usepackage{amsmath,amsfonts}
\usepackage{txfonts,wasysym,lscape,amssymb,mathrsfs,extarrows}
\usepackage[all,pdf]{xy}
\RequirePackage{color,graphicx}
\def\rev#1{{\color{black}{{#1}}}}

\newtheorem{theorem}{Theorem}[section]
\newtheorem{proposition}{Proposition}[section]
\newtheorem{lemma}{Lemma}[section]
\newtheorem{corollary}{Corollary}[section]
\newtheorem{definition}{Definition}[section]

\newtheorem{remark}{Remark}[section]

\title{Optimality Conditions for Constrained Minimax Optimization}

\author{Yu-Hong Dai\footnote{LSEC, ICMSEC, AMSS, Chinese Academy of Sciences, Beijing 100190, China.  {\sl Email}: dyh@lsec.cc.ac.cn.
This author was supported by the Natural Science Foundation of China (No. 11991020, 11631013, 11971372 and 11991021) and
Beijing Academy of Artificial Intelligence (BAAI).}
 \footnote{School of Mathematical Sciences, University of Chinese Academy of Sciences, Beijing 100049, China.}\quad and \quad Liwei Zhang \footnote{Corresponding author. School of Mathematical Sciences, Dalian University of Technology, Dalian 116024, China. {\sl Email}: lwzhang@dlut.edu.cn. This author was supported by the Natural Science Foundation of China (No. 11971089 and 11731013).}}

\date{}
\begin{document}

\maketitle

\begin{abstract}
Minimax optimization problems arises from both modern machine learning including generative adversarial
networks, adversarial training and multi-agent reinforcement learning, as well as from tradition research
areas such as saddle point problems, numerical partial differential equations and optimality conditions of
equality constrained optimization. For the unconstrained continuous nonconvex-nonconcave situation,
Jin, Netrapalli and Jordan (2019) carefully considered the very basic question: what is a proper
definition of local optima of a minimax optimization problem, and proposed a proper definition of
local optimality called local minimax. We shall extend the definition of local minimax point to constrained
nonconvex-nonconcave minimax optimization problems. By analyzing Jacobian uniqueness conditions for the lower-level maximization
problem and the strong regularity of Karush-Kuhn-Tucker conditions of the maximization problem, we provide
both necessary optimality conditions and sufficient optimality conditions for the local minimax points
of constrained minimax optimization problems.

\vskip 6 true pt \noindent \textbf{Key words}: constrained minimax optimization, value function, Jacobian uniqueness conditions, strong regularity, necessary optimality conditions, sufficient optimality conditions.
\vskip 12 true pt \noindent \textbf{AMS subject classification}: 90C30
\end{abstract}
\bigskip\noindent

\section{Introduction}
 \setcounter{equation}{0}
Minimax optimization problems arises from both modern machine learning including generative adversarial
networks, adversarial training and multi-agent reinforcement learning, as well as from tradition research
areas such as saddle point problems, numerical partial differential equations and optimality conditions of
equality constrained optimization.
Let $m,n,m_1,m_2,n_1$ and $n_2$ be positive integers, $f:\Re^n\times \Re^m \rightarrow \Re$, $h:\Re^n\times \Re^m \rightarrow \Re^{m_1}$, $g:\Re^n\times \Re^m \rightarrow \Re^{m_2}$, $H:\Re^n \rightarrow \Re^{n_1}$ and $G:\Re^n \rightarrow \Re^{n_2}$ be given functions. We are
interested in the constrained minimax optimization problem of the form
\begin{equation}\label{cminimax}
\min_{x \in \Phi}\max_{y \in Y(x)}f(x,y),
\end{equation}
where $f:\Re^n\times \Re^m \rightarrow \Re$, $\Phi\subset \Re^n$ is a feasible set of decision variable $x$ defined by
\begin{equation}\label{Phi}
\Phi=\{x \in \Re^n: H(x)=0, G(x)\leq 0\}
\end{equation}
and $Y: \Re^n \rightrightarrows \Re^m$ is a set-valued mapping defined by
\begin{equation}\label{Yx}
Y(x)=\{y \in \Re^m: h(x,y)=0, g(x,y)\leq 0\}.
\end{equation}

For the unconstrained continuous nonconvex-nonconcave situation, Jin, Netrapalli and Jordan (2019) \cite{Jin2019}
carefully considered the very basic question: what is a proper definition of local optima of a minimax optimization
problem, and proposed a proper definition of local optimality called local minimax. We shall extend this definition
of local minimax point for the constrained minimax optimization problem (\ref{cminimax}).
\begin{definition}\label{def:minimaxpoint}
A point $(x^*,y^*) \in \Re^n\times \Re^m$ is said to be a local minimax point of Problem (\ref{cminimax}) if  there exists $\delta_0>0$ and a function $\eta:(0,\delta_0] \rightarrow \Re_+$
 satisfying $\eta(\delta)\rightarrow 0$ as $\delta\rightarrow 0$, such that for any $\delta\in (0, \delta_0]$  and any $(x, y)\in [\textbf{B}_{\delta}(x^*)\cap \Phi] \times [Y(x^*) \cap \textbf{B}_{\delta}(y^*)]$, we have
 \begin{equation}\label{eq:definqs}
 f(x^*,y) \leq f(x^*,y^*)\leq \max\left\{f(x,z):z \in Y(x) \cap \textbf{B}_{\eta(\delta)}(y^*)\right\}.
 \end{equation}
\end{definition}
The minimax optimization problem is essentially a bi-level programming problem and  the local minimax point is \rev{closely} related to the so-called pessimistic solution of bi-level  programming problem, see \cite{Dempe2002}. There have been many results about  optimality conditions for bi-level  programming.  Dempe (1992)\cite{Dempe92} demonstrated necessary optimality conditions and the sufficient optimality conditions for the bi-level programming when the lower level problem is a convex optimization problem satisfying \rev{the Mangasarian-Fromovitz constraint qualification, the second-order sufficient optimality condition and the constant rank constraint qualification.}
Falk (1995)\cite{Falk95} discussed the optimality conditions when the lower level has a local unique solution and the upper level problem is unconstrained.
Ye and Zhu (1995,1997)\cite{YZhu1995} established necessary optimality conditions for bi-level programming based on the generalized gradient of value function.  Dempe et. al (2007)\cite{Dempe2007} derived necessary optimality conditions for bi-level programming when the solution set of the lower level problem satisfies some  calmness property. Dempe and Zemkoho (2013)\cite{Dempe2013} also developed necessary optimality conditions based on \rev{the} value function reformulation of bi-level programming and the assumption that the value function is \rev{locally convex}. Dempe et. al (2014)\cite{Dempe2014} derived a new type upper subdifferential necessary optimality conditions for the pessimistic version of bi-level
programming problem. Even recently, Mehlitz and Zemkoho (2019)\cite{MZemkoho2019} studied sufficient optimality conditions for bi-level programming.

Although many results about optimality conditions for bi-level programming are available, they are established based on different solution notations. For example, many works involve value function of the lower level problem, \rev{which is usually} defined as the global optimal value. \rev{This} restricts the application of the theoretical results. In this paper, we shall discuss optimality conditions for constrained minimax optimization whose solution is specified as the local minimax point given in Definition \ref{def:minimaxpoint}.

\section{Differential of the Value Function}\label{Sec2}
\setcounter{equation}{0}
\subsection{Under Jacobian uniqueness}
Let $(x^*,y^*)\in \Re^n \times \Re^m$ be a point and $f,h,g$ be twice continuously differentiable around $(x^*,y^*)$. For a point $x \in \Re^n$ around $x^*$, we use (${\rm P}_x$) to denote the following problem
\begin{equation}\label{eq:Px}
\begin{array}{cl}
\max_{z\in \Re^m} & f(x,z)\\[4pt]
{\rm s.t.\ \ \ \ \ }& h(x,z)=0,\\[4pt]
& g(x,z) \leq 0.
\end{array}
\end{equation}
The Lagrangian of Problem (${\rm P}_x$) is defined by
$$
{\cal L}(x;z,\mu,\lambda)=f(x,z)+\mu^Th(x,z)-\lambda^Tg(x,z).
$$
\begin{definition}\label{def:Jac}
Let $(\mu^*,\lambda^*) \in \Re^{m_1}\times \Re^{m_2}$ be a point. We say that Jacobian uniqueness conditions of Problem
(${\rm P}_{x^*}$) are satisfied at $(y^*,\mu^*,\lambda^*)$ if
\begin{itemize}
\item[(a)]The point $(y^*,\mu^*,\lambda^*)$ is a Karush-Kuhn-Tucker point of  Problem
(${\rm P}_{x^*}$); \rev{namely},
$$
\begin{array}{l}
\nabla_y{\cal L}(x^*;y^*,\mu^*,\lambda^*)=0,\\[3pt]
h(x^*,y^*)=0, \\[3pt]
 0 \leq \lambda^* \perp g(x^*,y^*)\leq 0.
\end{array}
$$
\item[(b)] \rev{The linear} independence constraint qualification holds at $y^*$; namely, the set of vectors
$$
\left\{\nabla_y h_1(x^*,y^*),\ldots,\nabla_y h_{m_1}(x^*,y^*)\right\} \cup \left\{\nabla_y g_i(x^*,y^*):i \in I_{x^*}(y^*)\right\}
$$
are linearly independent, where $I_{x^*}(y^*)=\left\{i: g_i(x^*,y^*)=0, i=1,\ldots,m_2\right\}$.
\item[(c)] \rev{The strict} complementarity condition holds at $y^*$ for $\lambda^*$; namely,
$$
\lambda^*_i-g_i(x^*,y^*)>0,\quad i=1, \ldots,m_2.
$$
\item[(d)] \rev{The second-order} sufficient optimality condition holds at $(y^*,\mu^*,\lambda^*)$,
$$
\langle \nabla_{yy}^2{\cal L}(x^*;y^*,\mu^*,\lambda^*)d_y,d_y\rangle<0 \quad \forall d_y \in {\cal C}_{x^*}(y^*),
$$
where ${\cal C}_{x^*}(y^*)$ is the critical cone of Problem
$({\rm P}_{x^*})$ at $y^*$,
$$
{\cal C}_{x^*}(y^*)=\left\{d_y \in \Re^m: {\cal J}_y h(x^*,y^*)d_y=0; \nabla_yg_i(x^*,y^*)d_y \leq 0, i \in I_{x^*}(y^*);
\nabla_yf(x^*,y^*)d_y\leq 0\right\}.
$$
\end{itemize}
\end{definition}
\begin{lemma}\label{lem:imp}
Let $(x^*,y^*) \in \Re^n \times \Re^m$ be a point around which $f,h,g$ are twice continuously differentiable. Let $(\mu^*,\lambda^*) \in \Re^{m_1}\times \Re^{m_2}$ such that Jacobian uniqueness conditions of Problem
$({\rm P}_{x^*})$ are satisfied at $(x^*,\mu^*,\lambda^*)$. Then there exist $\delta_0>0$ and $\varepsilon_0>0$, and a twice continuously differentiable mapping
$(y,\mu,\lambda):\textbf{B}_{\delta_0}(x^*) \rightarrow \textbf{B}_{\varepsilon_0}(y^*)\times \textbf{B}_{\varepsilon_0}(\mu^*)\times\textbf{B}_{\varepsilon_0}(\lambda^*)$ such that Jacobian uniqueness conditions of Problem $({\rm P}_x)$ are satisfied at $(y(x),\mu(x),\lambda(x))$ when $x \in \textbf{B}_{\delta_0}(x^*)$.
\end{lemma}
By introducing a set of auxiliary variables $w_1,\ldots, w_{m_2}$, we consider the following equality constrained optimization problem (${\rm \mathbb{P}}_x$).
\begin{equation}\label{AuxiP}
\begin{array}{cl}
\max & f(x,z)\\[4pt]
{\rm s.t.}& h(x,z)=0,\\[4pt]
& g(x,z) +w\circ w = 0,\\[4pt]
& z\in \Re^m,w \in \Re^{m_2},
\end{array}
\end{equation}
\rev{where $\circ$ is the Hadamard product.}

\begin{remark}\label{auxi-problem}
Let $w^*\in \Re^{m_2}$ with $w^*_i=\sqrt{-g_i(x^*,y^*)}$ for $i=1,\ldots,m_2$. With the help of Proposition 3.2 in \cite{B82}, we may prove that the point $(y^*,w^*,\mu^*,\lambda^*)$ satisfies Karush-Kuhn-Tucker conditions, \rev{the} linear independence constraint qualification as well as the second-order sufficient optimality condition for
 Problem (\ref{AuxiP}) hold at $(y^*,w^*,\mu^*,\lambda^*)$  if Jacobian uniqueness conditions of Problem
(${\rm P}_{x^*}$) are satisfied at $(x^*,\mu^*,\lambda^*)$.
\end{remark}
Define the optimal value function
\begin{equation}\label{eq:value-func}
\varphi (x)=f(x,y(x)),\quad x \in \textbf{B}_{\delta_0}(x^*),
\end{equation}
where $y(x)$ is defined by Lemma \ref{lem:imp}.
\begin{proposition}\label{prop:value-function}
\rev{If} the assumptions of Lemma \ref{lem:imp} are satisfied and $\varphi$ is defined by (\ref{eq:value-func}), then
\begin{equation}\label{phi-1deri}
\nabla_x \phi(x)=\nabla_x {\cal L}(x;y(x),\mu(x),\lambda(x))
\end{equation}
and
\begin{equation}\label{phi-2deri}
\begin{array}{ll}
\nabla^2 \phi(x)=& \nabla^2_{xx} {\cal L}(x;y(x),\mu(x),\lambda(x))\\[4pt]
&-\left[ \begin{array}{c}
\nabla^2_{yx}{\cal L}(x;y(x)\mu(x),\lambda(x))\\
 0\\
  {\cal J}_x h(x,y(x))\\
   {\cal J}_x g(x,y(x))
   \end{array}
   \right]^T K(x)^{-1}\left[ \begin{array}{c}
\nabla^2_{yx}{\cal L}(x;y(x),\mu(x),\lambda(x))\\
 0\\
  {\cal J}_x h(x,y(x))\\
   {\cal J}_x g(x,y(x))
   \end{array}
   \right],
   \end{array}
   \end{equation}
where
\begin{equation}\label{eq:Kx}
K(x)=
\left
[
\begin{array}{cccc}
\nabla^2_{yy}{\cal L}(x;y(x),\mu(x),\lambda(x)) & 0 & {\cal J}_y h(x,y(x))^T &
   {\cal J}_y g(x,y(x))^T\\[4pt]
   0  & -2 {\rm Diag}(\lambda(x)) & 0 & 2 {\rm Diag}\left(\sqrt{-g(x,y(x))}\right)\\[4pt]
   {\cal J}_y h(x,y(x)) & 0 &0 &0\\[4pt]
      {\cal J}_y g(x,y(x)) &2 {\rm Diag}\left(\sqrt{-g(x,y(x))}\right) & 0& 0
\end{array}
\right
].
\end{equation}
\end{proposition}
{\bf Proof.} Consider Problem (\ref{AuxiP}) and define its Lagrange function as
$$
\mathbb{L}(x;z,w,\mu,\lambda)=f(x,z)+\mu^Th(x,z)-\lambda^T(g(x,z)+w\circ w)={\cal L}(x;z,\mu,\lambda)-\lambda^Tw\circ w.
$$
From Remark \ref{auxi-problem}, we \rev{can} know that, when $x \in \textbf{B}_{\delta_0}(x^*)$, for $w(x)=\sqrt{-g(x,y(x))}$, the point
\rev{given by} $(y(x),w(x),\mu(x),\lambda(x))$ satisfies Kurash-Kuhn-Tucker conditions for Problem (${\rm \mathbb{P}}_x$), \rev{the} linear independence constraint qualification holds at $(y(x),w(x))$ and the second-order sufficient optimality condition holds at $(y(x),w(x),\mu(x),\lambda(x))$. Kurash-Kuhn-Tucker conditions for Problem (${\rm \mathbb{P}}_x$) at $(y(x),w(x),\mu(x),\lambda(x))$ can be expressed as
\begin{equation}\label{Auxi-KKT}
T(x(x),y(x),w(x),\mu(x),\lambda(x))=0,
\end{equation}
where
$$
T(x,y,w,\mu,\lambda)=
\rev{\left
[
\begin{array}{c}
\nabla_y\mathbb{L}(x;y,w,\mu,\lambda)\\[4pt]
\nabla_w\mathbb{L}(x;y,w,\mu,\lambda)\\[4pt]
h(x,y)\\[4pt]
g(x,y)+w\circ w
\end{array}
\right
]. }
$$
Differentiating both sides of (\ref{Auxi-KKT}) with respect to $x$ yields
\begin{equation}\label{eq:2side}
{\cal J}_{(y,w,\mu,\lambda)}T(x(x),y(x),w(x),\mu(x),\lambda(x))\left
[
\begin{array}{l}
{\cal J} y(x)\\
{\cal J} w(x)\\
{\cal J}\mu(x)\\
{\cal J}\lambda (x)
\end{array}\right
]+{\cal J}_{x}T(x(x),y(x),w(x),\mu(x),\lambda(x))=0.
\end{equation}
Noting that
$$
{\cal J}_{(y,w,\mu,\lambda)}T(x(x),y(x),w(x),\mu(x),\lambda(x))=K(x)
$$
and
$$
{\cal J}_{x}T(x(x),y(x),w(x),\mu(x),\lambda(x))=\left[ \begin{array}{c}
\nabla^2_{x,y}{\cal L}(x;y(x)\mu(x),\lambda(x))\\
 0\\
  {\cal J}_x h(x,y(x))\\
   {\cal J}_x g(x,y(x))
   \end{array}
   \right],
$$
we obtain (\ref{eq:2side}) that
\begin{equation}\label{eq:JD}
\left
[
\begin{array}{l}
{\cal J} y(x)\\
{\cal J} w(x)\\
{\cal J}\mu(x)\\
{\cal J}\lambda (x)
\end{array}\right]=-K(x)^{-1}\left[ \begin{array}{c}
\nabla^2_{x,y}{\cal L}(x;y(x),\mu(x),\lambda(x))\\
 0\\
  {\cal J}_x h(x,y(x))\\
   {\cal J}_x g(x,y(x))
   \end{array}
   \right].
\end{equation}
Noting that
$$
h(x,y(x))=0,\quad g(x,y(x))+w(x)\circ w(x)=0,
$$
we have
$$ \rev{
\begin{array}{rcl}
\varphi(x)&=&f(x,y(x))\\
 &=&f(x,y(x))+\mu(x)^Th(x,y(x))-\lambda(x)^T[g(x,y(x))+w(x)\circ w(x)]\\
 &=&\mathbb{L}(x;y(x),w(x),\mu(x),\lambda(x)).
\end{array}}
$$
Thus we get
$$\rev{
\begin{array}{rcl}
\nabla \varphi (x)& = & \nabla_x \mathbb{L}(x;y(x),w(x),\mu(x),\lambda(x))+{\cal J}y(x)^T\nabla_y\mathbb{L}(x;y(x),w(x),\mu(x),\lambda(x))\\[4pt]
& & +{\cal J}\mu(x)^Th(x,y(x))+{\cal J}\lambda(x)^T[g(x,y(x))+w(x)\circ w(x)]\\[4pt]
&= & \nabla_x \mathbb{L}(x;y(x),w(x),\mu(x),\lambda(x))\\
&= & \nabla_x {\cal L}(x;y(x),\mu(x),\lambda(x))
\end{array}}
$$
and
$$
\begin{array}{ll}
\nabla^2\varphi(x)&=\nabla^2_{xx}\mathbb{L}(x;y(x),w(x),\mu(x),\lambda(x))+{\cal J}_{y,w,\mu,\lambda}\nabla_x \mathbb{L}(x;y(x),w(x),\mu(x),\lambda(x)){\cal J}\left[
\begin{array}{l}
y(x)\\
w(x)\\
\mu(x)\\
\lambda(x)
\end{array}
\right]\\[4pt]
&=\nabla^2_{xx}\mathbb{L}(x;y(x),w(x),\mu(x),\lambda(x))\\[4pt]
& \,\,\,\,+\left[\nabla^2_{xy}\mathbb{L}(x;y(x),\mu(x),\lambda(x)) \quad 0\quad {\cal J}_xh(x,y(x))^T \quad {\cal J}_xg(x,y(x))\right]\left[
\begin{array}{l}
{\cal J}y(x)\\
{\cal J}w(x)\\
{\cal J}\mu(x)\\
{\cal J}\lambda(x)
\end{array}
\right]\\[4pt]
&=\nabla^2_{xx}{\cal L}(x;y(x),\mu(x),\lambda(x))\\[4pt]
& \,\,\,\,+\left[\nabla^2_{xy}\mathbb{L}(x;y(x),\mu(x),\lambda(x)) \quad 0\quad {\cal J}_xh(x,y(x))^T \quad {\cal J}_xg(x,y(x))\right]\left[
\begin{array}{l}
{\cal J}y(x)\\
{\cal J}w(x)\\
{\cal J}\mu(x)\\
{\cal J}\lambda(x)
\end{array}
\right].
\end{array}
$$
Combing this with (\ref{eq:JD}), we obtain (\ref{phi-1deri}) and (\ref{phi-2deri}). \hfill $\Box$

\subsection{Strong regularity of Karush-Kuhn-Tucker system}
We use $\Lambda_{x^*}(y^*)$  to denote the set of all $(\mu^*,\lambda^*) \in \Re^{m_1}\times \Re^{m_2}$ satisfying Karush-Kuhn-Tucker  conditions at $y^*$ for Problem
(${\rm P}_{x^*}$).
\begin{definition}\label{def:ssoc}
Let $(x^*,y^*) \in \Re^n \times \Re^m$ be a point at which $\Lambda_{x^*}(y^*) \ne \emptyset$. We say that the strong second-order sufficient optimality condition holds at $y^*$ for Problem
$({\rm P}_{x^*})$ if
$$
\sup_{(\mu,\lambda)\in \Lambda_{x^*}(y^*)}\langle \nabla_{yy}^2{\cal L}(x^*;y^*,\mu,\lambda)d_y,d_y\rangle<0 \quad \forall d_y \in {\rm aff\,}{\cal C}_{x^*}(y^*)\setminus \{0\},
$$
where ${\cal C}_{x^*}(y^*)$ is the critical cone of Problem
$({\rm P}_{x^*})$ at $y^*$.
\end{definition}
\begin{definition}\label{AssA}
Let $(x,y) \in \Re^n \times \Re^m$ be a point. We say that
{\bf Assumption A} holds at $y\in Y(x)$ for Problem $({\rm P}_x)$ if $\Lambda_x(y) \ne \emptyset$, \rev{the} linear independence constraint qualification and  the strong second-order sufficient optimality condition hold at $y$.
\end{definition}
\begin{lemma}\label{lem:LipImp}
Let $(x^*,y^*) \in \Re^n \times \Re^m$ be a point around which $f,h,g$ are twice continuously differentiable.
 Suppose that {\bf Assumption A}
  holds at $y^*$ for Problem
$({\rm P}_{x^*})$.
 Then there exist $\delta_0>0$ and $\varepsilon_0>0$, and a locally Lipschitz continuous mapping
$(y,\mu,\lambda):\textbf{B}_{\delta_0}(x^*) \rightarrow \textbf{B}_{\varepsilon_0}(y^*)\times \textbf{B}_{\varepsilon_0}(\mu^*)\times\textbf{B}_{\varepsilon_0}(\lambda^*)$ satisfying $(y(x^*),\mu(x^*),\lambda(x^*))=(y^*,\mu^*,\lambda^*)$ and
\begin{equation}\label{KKTs}
\begin{array}{l}
\nabla_y {\cal L}(x;y(x),\mu(x),\lambda(x))=0,\\[4pt]
h(x,y(x))=0,\\[4pt]
g(x,y(x))-\Pi_{\Re^{m_2}_-}(\lambda(x)+g(x,y(x)))=0
\end{array}
\end{equation}
for $x \in \textbf{B}_{\delta_0}(x^*)$. Moreover, {\bf Assumption A} holds at $y(x)$ for Problem $({\rm P}_x)$ when $x \in \textbf{B}_{\delta_0}(x^*)$.
\end{lemma}
{\bf Proof}. The first part of this lemma is from Robinson (1980) \cite{Robinson80}. For the second part, Karush-Kuhn-Tucker conditions of Problem (${\rm P}_x$) at $(y(x),\mu(x),\lambda(x))$ are from (\ref{KKTs}) and \rev{the} linear independence constraint qualification of $y(x)$ for Problem  (${\rm P}_x$) comes from the continuity of ${\cal J}_y h$ and ${\cal J}_yg$ and the linear independence constraint qualification  for (${\rm P}_{x^*}$) at $y^*$. We only need to prove the strong second-order sufficient optimality condition for Problem (${\rm P}_x$) at $(y(x),\mu(x),\lambda(x))$.
Let $\alpha_{x^*}=\{i: \lambda^*_i>0, i=1,\ldots, m_2\}$. Then the affine space of the critical cone  ${\cal C}_{x^*}(y^*)$ under {\bf Assumption A} is expressed as
$$
{\rm aff}\,{\cal C}_{x^*}(y^*)=\left\{d_y \in \Re^m: {\cal J}_y h(x^*,y^*)d_y=0; \nabla_y g_i(x^*,y^*)^Td_y=0, i \in \alpha_{x^*}\right\}.
$$
Noting that for $\alpha_{x}=\{i: \lambda_i (x) >0,i=1,\ldots, m_2\}$,
we have
$$
{\rm aff}\,{\cal C}_{x}(y(x))=\left\{d_y \in \Re^m: {\cal J}_y h(x,y(x))d_y=0; \nabla_y g_i(x,y(x))^Td_y=0, i \in \alpha_x\right\}.
$$
Since $\lambda (x)$ is continuous at $x^*$ and $\lambda^*_i>0$ for $\alpha_{x^*}$, we have $\alpha_x \supseteq \alpha_{x^*}$, and in turn,
\begin{equation}\label{eq:crx}
{\rm aff}\,{\cal C}_{x}(y(x))\subseteq \Gamma^*(x)=\left\{d_y \in \Re^m: {\cal J}_y h(x,y(x))d_y=0; \nabla_y g_i(x,y(x))^Td_y=0, i \in \alpha_{x^*}\right\}.
\end{equation}
It follows from the strong second-order sufficient condition for Problem
(${\rm P}_{x^*}$) at $y^*$ that \rev{the matrix} $\nabla^2_{yy}{\cal L}(x^*;y^*,\mu^*,\lambda^*)$ is negatively definite on $\Gamma^*(x^*)$. Thus we have for small $\delta_0>0$ and $x \in \textbf{B}_{\delta_0}(x^*)$ that  $\nabla^2_{yy}{\cal L}(x^*;y^*,\mu^*,\lambda^*)$ is  negatively definite on $\Gamma^*(x)$, \rev{which with (\ref{eq:crx}) implies that} $\nabla^2_{yy}{\cal L}(x^*;y^*,\mu^*,\lambda^*)$ is  negatively definite on ${\rm aff}\,{\cal C}_x(y(x))$. Finally, from the continuity of \rev{the matrix} $\nabla^2_{yy}{\cal L}(x;y(x),\mu(x),\lambda(x))$ with respect to $x$ around $x^*$, we obtain that $\nabla^2_{yy}{\cal L}(x;y(x),\mu(x),\lambda(x))$ is  negatively definite over ${\rm aff}\,{\cal C}_x(y(x))$, \rev{indicating}  the strong second-order sufficient optimality condition holds at $y(x)$ for (${\rm P}_{x}$)
 when $x$ is around $x^*$. \hfill $\Box$
\begin{remark}\label{rem:isoM}
Let {\bf Assumption A} be satisfied. Then from Lemma \ref{lem:LipImp}, for $x \in \textbf{B}_{\delta_0}(x^*)$, $y(x)$ is a unique local maximizer of Problem $({\rm P}_{x})$.
\end{remark}
For a linear operator $W: \Re^{m_2} \rightarrow \Re^{m_2}$, define
\begin{equation}\label{eq:amatrix}
{\cal A}(x,W)=\left
[
\begin{array}{ccc}
\nabla^2_{yy}{\cal L}(x; y(x),\mu(x),\lambda(x)) & {\cal J}_yh(x,y(x))^T & {\cal J}_yg(x,y(x))^T\\[4pt]
 {\cal J}_yh(x,y(x)) & 0 & 0 \\[4pt]
  (I-W){\cal J}_yg(x,y(x)) &0 & W
  \end{array}
  \right
  ],
\end{equation}
 a set-valued mapping $\mathbb{A}_B:\Re^n \rightrightarrows \Re^{m_2 \times m_2}$ by
 \begin{equation}\label{def:Aset}
 \mathbb{A}_B(x)=\left\{{\cal A}(x,W):W \in \partial_B \Pi_{\Re^{m_2}_-}(\lambda(x)+g(x,y(x)))\right\},
 \end{equation}
 and a set-valued mapping $\mathbb{A}_C:\Re^n \rightrightarrows \Re^{m_2 \times m_2}$ by
 \begin{equation}\label{def:Cset}
 \mathbb{A}_C(x)=\left\{{\cal A}(x,W):W \in \partial\Pi_{\Re^{m_2}_-}(\lambda(x)+g(x,y(x)))\right\}.
 \end{equation}
 \begin{proposition}\label{prop:AnonS}
 Let $(x^*,y^*) \in \Re^n \times \Re^m$ be a point around which $f,h,g$ are twice continuously differentiable.
 Suppose that {\bf Assumption A}
  holds at $y^*$ for Problem
$({\rm P}_{x^*})$. Then every element in $\mathbb{A}_C(x^*)$ is nonsingular.
 \end{proposition}
 {\bf Proof}.  Let $V\in \mathbb{A}_B(x^*)$. Then there exists an element $W \in \partial_B \Pi_{\Re^{m_2}_-}(\lambda^*+g(x^*,y^*))$ such that
 $$
V=\left
[
\begin{array}{ccc}
\nabla^2_{yy}{\cal L}(x^*; y^*,\mu^*,\lambda^*) & {\cal J}_yh(x^*,y^*)^T & {\cal J}_yg(x^*,y^*))^T\\[4pt]
 {\cal J}_yh(x^*,y^*) & 0 & 0 \\[4pt]
  (I-W){\cal J}_yg(x^*y^*) &0 & W
  \end{array}
  \right
  ].
 $$
 Define
 $$
 \alpha=\{i: g_i(x^*,y^*)=0,\lambda_i>0\},\ \beta=\{i: g_i(x^*,y^*)=0,\lambda_i=0\},\ \gamma=\{i: g_i(x^*,y^*)<0,\ \lambda_i=0\}.
 $$
 Hence $W$ can be expressed as
 $$
 W={\rm Diag}(w_1,\ldots,w_{m_2})
 $$
 with
 $$
 w_i \left
 \{
 \begin{array}{ll}
 =0  & i \in \alpha,\\[2pt]
 \in [0,1] & i \in \beta,\\[2pt]
 =1 &  i \in  \gamma.
 \end{array}
 \right.
 $$
 Then $V$ may be expressed as
 $$
 V=\left
[
\begin{array}{ccccc}
\nabla^2_{yy}{\cal L}(x^*; y^*,\mu^*,\lambda^*) & {\cal J}_yh(x^*,y^*)^T & {\cal J}_yg_{ \alpha}(x^*,y^*)^T& {\cal J}_yg_{\beta}(x^*,y^*)^T & {\cal J}_yg_{ \gamma}(x^*,y^*)^T\\[4pt]
 {\cal J}_yh(x^*,y^*) & 0 & 0 &0 &0\\[4pt]
  {\cal J}_yg_{\alpha}(x^*y^*) &0 & 0 & W_{\beta} &0\\[4pt]
  (I_{\beta}-W_{\beta}){\cal J}_yg_{\beta}(x^*y^*) &0 & 0 &0 &0\\[4pt]
  0 & 0 &0 & 0& I_{|\gamma|}
  \end{array}
  \right
  ].
  $$
 For $\xi_1 \in \Re^m$, $\xi_2 \in \Re^{m_1}$,$\xi_3\in \Re^{|\alpha|}$, $\xi_4\in \Re^{|\beta|}$ and $\xi_5 \in \Re^{| \gamma|}$, consider
 $$
 V\left[
 \begin{array}{c}
 \xi_1\\
 \xi_2\\
 \xi_3\\
 \xi_4\\
 \xi_5
 \end{array}
 \right]=0.
 $$
 or equivalently
 $$
\begin{array}{rr}
\nabla^2_{yy}{\cal L}(x^*; y^*,\mu^*,\lambda^*) \xi_1+ {\cal J}_yh(x^*,y^*)^T\xi_2+ {\cal J}_yg_{\alpha}(x^*,y^*)^T\xi_3\,\\[2pt]
+ {\cal J}_yg_{\beta}(x^*,y^*)^T\xi_4+ {\cal J}_yg_{\gamma}(x^*,y^*)^T\xi_5& =0,\\[2pt]
 {\cal J}_yh(x^*,y^*)\xi_1&=0, \\[2pt]
  {\cal J}_yg_{\alpha}(x^*y^*)\xi_1&=0,\\[2pt]
  (I_{|\beta|}-W_{\beta}){\cal J}_yg_{\beta}(x^*y^*)\xi_1+W_{\beta}\xi_4&=0,\\[2pt]
  \xi_5&=0,
   \end{array}
 $$
 which implies $\xi_5=0$ and
 \begin{equation}\label{eq:help1}
\begin{array}{r}
\nabla^2_{yy}{\cal L}(x^*; y^*,\mu^*,\lambda^*) \xi_1+ {\cal J}_yh(x^*,y^*)^T\xi_2+ {\cal J}_yg_{\alpha}(x^*,y^*)^T\xi_3 + {\cal J}_yg_{\beta}(x^*,y^*)^T\xi_4 =0,\\[2pt]
 {\cal J}_yh(x^*,y^*)\xi_1=0, \\[2pt]
  {\cal J}_yg_{ \alpha}(x^*y^*)\xi_1 =0,\\[2pt]
  (I_{|\beta|}-W_{\beta}){\cal J}_yg_{\beta}(x^*y^*)\xi_1+W_{\beta}\xi_4=0.
     \end{array}
 \end{equation}
 Without loss of generality, we assume that $w_i \in (0, 1]$ so that $W_{\beta}$ is nonsingular. From the fourth row  of (\ref{eq:help1}),we have
 \begin{equation}\label{eq:hp2}
 \xi_4=-W_{\beta}^{-1}(I_{|\beta|}-W_{\beta}){\cal J}_y g_{\beta}(x^*,y^*)\xi_1.
 \end{equation}
 From the second row and the third row, we have
 $$
 {\cal J}_y h(x^*,y^*) \xi_1=0, {\cal J}_yg_{\alpha}(x^*,y^*) \xi_1=0,
 $$
 which implies $\xi_1 \in {\rm aff}\, {\cal C}_{x^*}(y^*)$. Substituting $\xi_4$ of (\ref{eq:hp2})  in (\ref{eq:help1}) and premultiplying $\xi_1^T$ to both sides of the first row of (\ref{eq:help1}), we obtain
 \begin{equation}\label{eq:hp3}
\xi_1^T \nabla^2_{yy}{\cal L}(x^*; y^*,\mu^*,\lambda^*) \xi_1  -\xi_1^T{\cal J}_y g_{\beta}(x^*,y^*)^TW_{\beta}^{-1}(I_{|\beta|}-W_{\beta}){\cal J}_y g_{\beta}(x^*,y^*)\xi_1=0.
 \end{equation}
 Since $TW_{\beta}^{-1}(I_{|\beta|}-W_{\beta})$ is diagonal with $W_{ii}=(1-w_i)/w_i>0$, we have that this matrix is positively definite, and
 $$
  -\xi_1^T{\cal J}_y g_{\beta}(x^*,y^*)^TW_{\beta}^{-1}(I_{|\beta|}-W_{\beta}){\cal J}_y g_{\beta}(x^*,y^*)\xi_1\leq0.
 $$
 Therefore we obtain $\xi_1=0$ from  $\xi_1 \in {\rm aff}\, {\cal C}_{x^*}(y^*)$ and the strong second-order sufficient optimality condition. Obviously we have $\xi_4=0$. Substitute $\xi_1=0$ and $\xi_4=0$ to the first row of (\ref{eq:help1}) and using linear independence constraint qualification, we get $\xi_2=0$ and $\xi_3=0$. Therefore, matrix $V$ is nonsingular. The proof is complete. \hfill $\Box$
\begin{corollary}\label{coro:AnonS}
 Let $(x^*,y^*) \in \Re^n \times \Re^m$ be a point around which $f,h,g$ are twice continuously differentiable.
 Suppose that {\bf Assumption A}
  holds at $y^*$ for Problem
$({\rm P}_{x^*})$. Then every element in $\mathbb{A}_B(x^*)$ is nonsingular.
 \end{corollary}
 {\bf Proof}.  Let $V\in \mathbb{A}_B(x^*)$. Then there exists an element $W \in \partial_B \Pi_{\Re^{m_2}_-}(\lambda^*+g(x^*,y^*))$ such that
 $$
V=\left
[
\begin{array}{ccc}
\nabla^2_{yy}{\cal L}(x^*; y^*,\mu^*,\lambda^*) & {\cal J}_yh(x^*,y^*)^T & {\cal J}_yg(x^*,y^*))^T\\[4pt]
 {\cal J}_yh(x^*,y^*) & 0 & 0 \\[4pt]
  (I-W){\cal J}_yg(x^*y^*) &0 & W
  \end{array}
  \right
  ].
 $$
 Define
 $$
 \alpha=\{i: g_i(x^*,y^*)=0,\lambda_i>0\}, \beta=\{i: g_i(x^*,y^*)=0,\lambda_i=0\},\gamma=\{i: g_i(x^*,y^*)<0,\lambda_i=0\}.
 $$
 There exists a partition of $\beta$, say $(\beta_1,\beta_2)$, namely $\beta_1\cup\beta_2 =\beta$ and $\beta_1\cap \beta_2=\emptyset$, such that
 $$
 W={\rm Diag}(w_1,\ldots,w_{m_2})
 $$
 with
 $$
 w_i =\left
 \{
 \begin{array}{ll}
 0 & i \in \tilde \alpha ,\\[2pt]
 1 &  i \in \tilde \gamma,
 \end{array}
 \right.
 $$
 where $\tilde \alpha =\alpha \cup \beta_1$ and $\tilde \gamma=\beta_2 \cup \gamma$.  Then $V$ may be expressed as
 $$
 V=\left
[
\begin{array}{cccc}
\nabla^2_{yy}{\cal L}(x^*; y^*,\mu^*,\lambda^*) & {\cal J}_yh(x^*,y^*)^T & {\cal J}_yg_{\tilde \alpha}(x^*,y^*)^T& {\cal J}_yg_{\tilde \gamma}(x^*,y^*)^T\\[4pt]
 {\cal J}_yh(x^*,y^*) & 0 & 0 &0 \\[4pt]
  {\cal J}_yg_{\tilde \alpha}(x^*y^*) &0 & 0 &0\\[4pt]
  0 & 0 &0 & I_{|\tilde \gamma|}
  \end{array}
  \right
  ].
  $$
  The nonsingulairity of $V$ can be proved in a similar way as in the proof of Proposition \ref{prop:AnonS}.\hfill $\Box$
 \begin{proposition}\label{prop:Anonx}
 Let $(x^*,y^*) \in \Re^n \times \Re^m$ be a point around which $f,h,g$ are twice continuously differentiable.
 Suppose that {\bf Assumption A}
  holds at $y^*$ for Problem
$({\rm P}_{x^*})$.
 Let  $\delta_0>0$ be given in Lemma \ref{lem:LipImp}. Then the set-valued mapping
 $\mathbb{A}_B(x)$ $(\mathbb{A}_C(x))$ is upper semicontinuous at $x^*$, and
 for small $\delta \in (0,\delta_0)$, every element in $\mathbb{A}_B(x)$ $(\mathbb{A}_C(x))$ is nonsingular when $x \in \textbf{B}(x^*,\delta)$.
 \end{proposition}
 \begin{proposition}\label{prop:dirD}
  Let $(x^*,y^*) \in \Re^n \times \Re^m$ be a point around which $f,h,g$ are twice continuously differentiable.
 Suppose that {\bf Assumption A}
  holds at $y^*$ for Problem
$({\rm P}_{x^*})$. Let  $\delta_0>0$, $\varepsilon_0>0$ and $(y(\cdot),\mu(\cdot),\lambda(\cdot))$ be given in Lemma \ref{lem:LipImp}. Then for $x \in \textbf{B}_{\delta_0}(x^*)$,
\begin{itemize}
\item[(i)]The directional derivative of $(y(\cdot),\mu(\cdot),\lambda(\cdot))$ at $x$  satisfies
\begin{equation}\label{eq:DD}
\left
(
\begin{array}{c}
y'(x;d_x)\\
\mu'(x;d_x)\\
\lambda'(x;d_x)
\end{array}
\right
)\in \left\{-
{\cal A}(x,W)^{-1}\left
(
\begin{array}{c}
\nabla^2_{yx}{\cal L}(x;y(x),\mu(x),\lambda(x))d_x\\
{\cal J}_x h(x,y(x))d_x\\
(I-W){\cal J}_xg(x,y(x))d_x
\end{array}
\right
): W \in \partial_B \Pi_{\Re^{m_2}_-}(\lambda(x)+g(x,y(x)))
\right \}.
\end{equation}
\item[(ii)]The B-subdifferential of $(y(\cdot),\mu(\cdot),\lambda(\cdot))$ at $x$  satisfies
\begin{equation}\label{eq:BD}
\partial_B \left
(
\begin{array}{c}
y\\
\mu\\
\lambda
\end{array}
\right
)(x)\in \left\{-
{\cal A}(x,W)^{-1}\left
(
\begin{array}{c}
\nabla^2_{yx}{\cal L}(x;y(x),\mu(x),\lambda(x))\\
{\cal J}_x h(x,y(x))\\
(I-W){\cal J}_xg(x,y(x))
\end{array}
\right
): W \in \partial_B \Pi_{\Re^{m_2}_-}(\lambda(x)+g(x,y(x)))
\right \}.
\end{equation}
\item[(iii)]Clarke generalized Jacobian of $(y(\cdot),\mu(\cdot),\lambda(\cdot))$ at $x$  satisfies
\begin{equation}\label{eq:CD}
\partial \left
(
\begin{array}{c}
y\\
\mu\\
\lambda
\end{array}
\right
)(x)\in \left\{
{\cal A}(x,W)^{-1}\left
(
\begin{array}{c}
\nabla^2_{yx}{\cal L}(x;y(x),\mu(x),\lambda(x))\\
{\cal J}_x h(x,y(x))\\
(I-W){\cal J}_xg(x,y(x))
\end{array}
\right
): W \in \partial\Pi_{\Re^{m_2}_-}(\lambda(x)+g(x,y(x)))
\right \}.
\end{equation}
\end{itemize}
 \end{proposition}
 {\bf Proof}.  It follows from (\ref{KKTs}), for $d_x \in \Re^n$ that
 \begin{equation}\label{KKTsDer}
\begin{array}{rr}
\nabla^2_{yx} {\cal L}(x;y(x),\mu(x),\lambda(x))d_x
+\nabla^2_{yy} {\cal L}(x;y(x),\mu(x),\lambda(x))y'(x;d_x)&\\
\quad +{\cal J}_yh(x,y(x))\mu'(x;d_x)-
{\cal J}_yg(x,y(x))^T\lambda'(x;d_x)&=0,\\[4pt]
{\cal J}_xh(x,y(x))d_x+{\cal J}_yh(x,y(x))y'(x;d_x)&=0,\\[4pt]
{\cal J}_xg(x,y(x))d_x+{\cal J}_yg(x,y(x))y'(x;d_x)-\Pi'_{\Re^{m_2}_-}(\lambda(x)+g(x,y(x));\lambda'(x;d_x)&\\
 +{\cal J}_xg(x,y(x))d_x+{\cal J}_yg(x,y(x))y'(x;d_x))&=0.
\end{array}
\end{equation}
 Since $\Pi_{\Re^{m_2}_-}$ is semismooth everywhere, we have from \cite{Qi93a} that there exists an matrix
 $\widetilde W \in \partial_B\Pi_{\Re^{m_2}_-}(\lambda(x)+g(x,y(x)))$ such that
 \begin{equation}\label{eq:derW}
 \begin{array}{r}
 \Pi'_{\Re^{m_2}_-}(\lambda(x)+g(x,y(x));\lambda'(x;d_x)+{\cal J}_xg(x,y(x))d_x+{\cal J}_yg(x,y(x))y'(x;d_x))\\[6pt]
 \quad \,=\widetilde{W}[\lambda'(x;d_x)+{\cal J}_xg(x,y(x))d_x+{\cal J}_yg(x,y(x))y'(x;d_x)].
 \end{array}
 \end{equation}
 Substituting (\ref{eq:derW}) into (\ref{KKTsDer}), we may rewrite (\ref{KKTsDer}) as
 $$
{\cal A}(x,\widetilde W) \left
(
\begin{array}{c}
y'(x;d_x)\\
\mu'(x;d_x)\\
\lambda'(x;d_x)
\end{array}
\right
)=-
\left
(
\begin{array}{c}
\nabla^2_{yx}{\cal L}(x;y(x),\mu(x),\lambda(x))d_x\\
{\cal J}_x h(x,y(x))d_x\\
(I-\widetilde W){\cal J}_xg(x,y(x))d_x
\end{array}
\right
),
$$
which implies (\ref{eq:DD}).

 We use $(y(x);\mu(x);\lambda (x))$ to denote $(y(x)^T,\mu(x)^T\lambda (x)^T)^T$. Let ${\cal D}_1(\delta_0)$ be the set of all differentiable  points of  $(y(x);\mu(x);\lambda (x))$ in $\textbf{B}_{\delta_0}(x^*)$ and ${\cal D}_2(\delta_0)$ be the set of all differentiable  points such that $\lambda(x)-g(x,(y(x))>0$ in $\textbf{B}_{\delta_0}(x^*)$ (namely; the points at which $\Pi_{\Re^{m_2}_-}$ is differentiable), and define ${\cal D}(\delta_0)={\cal D}_1(\delta_0)\cap {\cal D}_2(\delta_0)$.
For every $V\in \partial_B (y(x);\mu(x);\lambda (x))$, there exists a sequence $x^k \rightarrow x$ with $x^k \in {\cal D}(\delta_0)$ such that
$$
{\cal J}_x (y(x^k);\mu(x^k);\lambda (x^k))\rightarrow V.
$$
It follows from (\ref{KKTs}) and $x^k \in \in {\cal D}(\delta_0)$   that
 $$
\begin{array}{rr}
\nabla^2_{yx} {\cal L}(x^k;y(x^k),\mu(x^k),\lambda(x^k))
+\nabla^2_{yy} {\cal L}(x^k;y(x^k),\mu(x^k),\lambda(x^k)){\cal J}y(x^k)&\\[4pt]
\quad +{\cal J}_yh(x^k,y(x^k)){\cal J}\mu(x^k)-
{\cal J}_yg(x^k,y(x^k))^T{\cal J}\lambda(x^k)&=0,\\[4pt]
{\cal J}_xh(x^k,y(x^k))+{\cal J}_yh(x^k,y(x^k)){\cal J}y(x^k)&=0,\\[4pt]
{\cal J}_xg(x^k,y(x^k))+{\cal J}_yg(x^k,y(x^k)){\cal J}y(x^k)-{\cal J}\Pi_{\Re^{m_2}_-}(\lambda(x^k)+g(x^k,y(x^k)))({\cal J}\lambda(x^k)&\\[4pt]
+{\cal J}_xg(x^k,y(x^k))+{\cal J}_yg(x^k,y(x^k)){\cal J}y(x^k))&=0
\end{array}
$$
or, equivalently,
 \begin{equation}\label{KKTsgrdk}
{\cal J} \left
(
\begin{array}{c}
y\\
\mu\\
\lambda
\end{array}
\right
)(x^k)=-{\cal A}(x^k, W^k)^{-1}
\left
(
\begin{array}{c}
\nabla^2_{yx}{\cal L}(x^k;y(x^k),\mu(x^k),\lambda(x^k))\\
{\cal J}_x h(x^k,y(x^k))\\
(I- W^k){\cal J}_xg(x^k,y(x^k))
\end{array}
\right
)
\end{equation}
 with
 $$
 W^k={\cal J}\Pi_{\Re^{m_2}_-}(\lambda(x^k)+g(x^k,y(x^k))).
 $$
 Let $W=\lim_{\rev{k\rightarrow \infty}} W^k$ (or assume that $W$ is an limit operator of $W^k$). Then $W \in \partial_B \Pi_{\Re^{m_2}_-}(\lambda(x)+g(x,y(x)))$.  Taking the limit in both sides of (\ref{KKTsgrdk}) as $k \rightarrow \infty$, we obtain the result in (\ref{eq:BD}). From the definition of Clarke generalized Jacobian, we obtain (\ref{eq:CD}) from (\ref{eq:BD}). \hfill $\Box$\\
 Define
 \begin{equation}\label{eq:Hxw}
 H(x,W)={\cal A}(x,W)^{-1}\left
(
\begin{array}{c}
\nabla^2_{yx}{\cal L}(x;y(x),\mu(x),\lambda(x))\\
{\cal J}_x h(x,y(x))\\
(I-W){\cal J}_xg(x,y(x))
\end{array}
\right
).
 \end{equation}
 Then we obtain the directional derivative, B-subdifferential and Clarke generalized subdifferential of $\phi$ at $x$ by the following corollary.
\begin{corollary}\label{phi-d}
Let $(x^*,y^*) \in \Re^n \times \Re^m$ be a point around which $f,h,g$ are twice continuously differentiable.
 Suppose that {\bf Assumption A}
  holds at $y^*$ for Problem
$({\rm P}_{x^*})$. Let  $\delta_0>0$, $\varepsilon_0>0$ and $(y(\cdot),\mu(\cdot),\lambda(\cdot))$ be given in Lemma \ref{lem:LipImp}. Then $\varphi$ is locally Lipschitz continuous in $\textbf{B}_{\delta_0}(x^*)$ and for $x \in \textbf{B}_{\delta_0}(x^*)$,
\begin{itemize}
\item[$(i)$]The directional derivative of $\varphi$ at $x$  satisfies
\begin{equation}\label{eq:objD}
\begin{array}{l}
\varphi'(x;d_x)\in \nabla_x{\cal L}(x;y(x),\mu(x),\lambda(x))d_x\\
\quad  \,-
 \left\{\nabla_{y,\mu,\lambda}{\cal L}(x;y(x),\mu(x),\lambda(x))^TH(x,W)d_x
: W \in \partial_B \Pi_{\Re^{m_2}_-}(\lambda(x)+g(x,y(x)))
\right \}.
\end{array}
\end{equation}
\item[$(ii)$]The  B-subdifferential of $\varphi$ at $x$  satisfies
\begin{equation}\label{eq:objBD}
\begin{array}{l}
\partial_B\varphi(x) \in \nabla_x{\cal L}(x;y(x),\mu(x),\lambda(x))\\
\quad  \,-\left\{
 H(x,W)^T
: W \in \partial_B \Pi_{\Re^{m_2}_-}(\lambda(x)+g(x,y(x)))
\right \}\nabla_{y,\mu,\lambda}{\cal L}(x;y(x),\mu(x),\lambda(x)).
\end{array}
\end{equation}
\item[$(ii)$]The  Clarke generalized subdifferential of $\varphi$ at $x$  satisfies
\begin{equation}\label{eq:objCD}
\begin{array}{l}
\partial\varphi(x) \in \nabla_x{\cal L}(x;y(x),\mu(x),\lambda(x))\\
\quad  \,-\left\{
 H(x,W)^T
: W \in \partial \Pi_{\Re^{m_2}_-}(\lambda(x)+g(x,y(x)))
\right \}\nabla_{y,\mu,\lambda}{\cal L}(x;y(x),\mu(x),\lambda(x)).
\end{array}
\end{equation}
\end{itemize}
\end{corollary}

\section{Optimality Conditions}
\setcounter{equation}{0}
\quad Suppose that $\varphi (x)$ is defined by (\ref{eq:value-func}). Then the constrained minimax problem (\ref{cminimax}) is locally reduced to
\begin{equation}\label{prob:RP}
\begin{array}{cl}
\min & \varphi (x)=f(x,y(x))\\[3pt]
{\rm s.t.} & x \in \Phi \cap \textbf{B}_{\delta_0}(x^*),
\end{array}
\end{equation}
where $y(x)$ is a local minimizer of (${\rm P}_x$) around $y^*$ and  $\Phi$ is defined by (\ref{Phi}).

For $x^* \in \Phi$, the Mangasarian-Fromovitz constraint qualification is said to hold at $x^*$ the constraint set $\Phi$ if
\begin{itemize}
\item[(a)] The set of vectors $\nabla H_j(x^*), j=1,\ldots, n_1$ are linearly independent;
\item[(b)] There exists a vector $\bar d \in \Re^n$ such that
$$
\nabla H_j(x^*)^T\bar d=0,j=1,\ldots,n_1, \nabla G_i(x^*)^T\bar d<0, i \in I(x^*),
$$
where $I(x^*)=\{i: G_i(x^*)=0, i=1,\ldots, n_2\}$.
\end{itemize}
Define the critical cone of Problem (\ref{prob:RP}) at $x^*$ is defined by
\begin{equation}\label{eq:critC}
{\cal C}(x^*)=\{d_x \in \Re^n: {\cal J}H(x^*)d_x=0; \nabla G_i(x^*)^Td_x \leq 0, i \in I(x^*); \varphi'(x^*;d_x) \leq0\}.
\end{equation}

We now derive necessary optimality conditions and second-order sufficient optimality conditions for Problem (\ref{cminimax}) under Jacobian uniqueness conditions for (${\rm P}_{x^*}$).  In this case, the critical cone ${\cal C}(x^*)$ can be expressed as
\begin{equation}\label{eq:CritiCS}
{\cal C}(x^*)=\{d_x \in \Re^n: {\cal J}H(x^*)d_x=0; \nabla G_i(x^*)^Td_x \leq 0, i \in I(x^*); \nabla_x {\cal L}(x^*;y^*,\mu^*,\lambda^*)^Td_x\leq 0\}.
\end{equation}
\begin{theorem}\label{th:nc}
(Necessary Optimality Conditions) Let $(x^*,y^*) \in \Re^n \times \Re^m$ be a point around which $f,h,g$ are twice continuously differentiable and $H$, $G$ are twice continuously differentiable around $x^*$. Let $(x^*,y^*)$  be a local minimax point of Problem (\ref{cminimax}). Assume  that the linear independence constraint qualification holds at $y^*$ for constraint set $Y(x^*)$. Then there exists a unique vector $(\mu^*,\lambda^*) \in \Re^{m_1}\times\Re^{m_2}$ such that
\begin{equation}\label{KKT-Px}
\begin{array}{l}
\nabla_y {\cal L}(x^*;y^*,\mu^*,\lambda^*)=0,\\[3pt]
h(x^*,y^*)=0,\\[3pt]
0\geq \lambda^* \perp g(x^*,y^*) \leq 0.
\end{array}
\end{equation}
For any $d_y \in {\cal C}_{x^*}(y^*)$, we have that
\begin{equation}\label{second-N}
\langle \nabla^2_{yy}{\cal L}(x^*;y^*,\mu^*,\lambda^*)d_y, d_y \rangle \leq 0.
\end{equation}
 Assuming Problem $({\rm P}_{x^*})$  satisfies  Jacobian uniqueness conditions at $(y^*,\mu^*,\lambda^*)$ and
  the Mangasarian-Fromovitz constraint qualification holds at $x^*$ for the constraint set $\Phi$,
 there exists $(u^*,v^*) \in \Re^{n_1}\times \Re^{n_2}$ such that
 \begin{equation}\label{KKT-P}
\begin{array}{l}
\nabla_x {\cal L}(x^*;y^*,\mu^*,\lambda^*)+{\cal J}H(x^*)^Tu^*+{\cal J}G(x^*)^Tv^*=0,\\[3pt]
H(x^*)=0,\\[3pt]
0\leq v^* \perp G(x^*) \leq 0.
\end{array}
\end{equation}
The set of all $(u^*,v^*)$ satisfying (\ref{KKT-P}), denoted by $\Lambda (x^*)$, is nonempty compact convex set.
Furthermore, for every $d_x \in {\cal C}(x^*)$ where ${\cal C}(x^*)$ is defined by (\ref{eq:CritiCS}),
\begin{equation}\label{eq:secNCs}
\begin{array}{l}
\displaystyle \max_{(u,v) \in \Lambda(x^*)} \left\{ \left\langle \left[\displaystyle \sum_{j=1}^{n_1}u_i\nabla^2_{xx}H_j(x^*)+\displaystyle \sum_{i=1}^{n_2} v_i\nabla^2_{xx}G_i(x^*)\right]d_x,d_x \right\rangle\right\}\\[16pt]
 \quad \quad \quad \, +\left\langle \left[ \nabla^2_{xx}{\cal L}(x^*;y^*,\mu^*,\lambda^*)- N(x^*)^TK(x^*)^{-1}N(x^*)\right]d_x,d_x \right\rangle\geq0,
 \end{array}
\end{equation}
where $K(x)$ is defined by (\ref{eq:Kx}) and $N(x)$ is defined by
\begin{equation}\label{Nx}
N(x)=\left[ \begin{array}{c}
\nabla^2_{x,y}{\cal L}(x;y(x)\mu(x),\lambda(x))\\
 0\\
  {\cal J}_x h(x,y(x))\\
   {\cal J}_x g(x,y(x))
   \end{array}
   \right].
\end{equation}
\end{theorem}
{\bf Proof}. Since $y^*$ is a local minimizer of (${\rm P}_{x^*}$) and the linear independence constraint qualification holds at $y^*$ for the constraint set
$$
Y(x^*)=\{y \in \Re^m: h(x^*,y)=0,\,\,g(x^*,y) \leq 0\},
$$
we may obtain the first-order and second-order necessary optimality conditions (\ref{KKT-Px}) and (\ref{second-N}) from \cite{NW99}.  Noting that $x^*$ is a local minimizer of the following problem
$$
\begin{array}{cl}
\min & \varphi (x)=f(x,y(x))\\[3pt]
{\rm s.t.} & x \in \Phi \cap \textbf{B}_{\delta_0}(x^*).
\end{array}
$$
The Lagrange function of the above problem is
$$
L(x,u,v)=\varphi(x)+u^TH(x)+v^TG(x).
$$
It follows from \cite{BS00} that there exist $u^*$ and $v^*$ such that
\begin{equation}\label{eq:KKTP}
\begin{array}{l}
\nabla_x L(x^*,u^*,v^*)=0,\\[4pt]
H(x^*)=0,\\[4pt]
0\leq u^* \perp G(x^*)\leq 0.
\end{array}
\end{equation}
Since the Mangasarian-Fromovitz constraint qualification holds at $x^*$ for the constraint set $\Phi$, the set of all vectors $(u^*,v^*)$ satisfying (\ref{eq:KKTP}) is a nonempty compact convex set. From the formula for $\nabla \varphi (x)$ in (\ref{phi-1deri}), we obtain (\ref{KKT-P}) from (\ref{eq:KKTP}) and $\Lambda (x^*)$  is nonempty compact convex.  It also follows from \cite{BS00} that the second-order necessary optimality conditions at $x^*$ can be expressed as
\begin{equation}\label{eq:SONC}
\displaystyle \max_{(u,v) \in \Lambda(x^*)} \left\{ \left\langle \left[\nabla^2 \varphi (x^*)+\displaystyle \sum_{j=1}^{n_1}u_i\nabla^2_{xx}H_j(x^*)+\displaystyle \sum_{i=1}^{n_2} v_i\nabla^2_{xx}G_i(x^*)\right]d_x,d_x \right\rangle\right\}\geq 0, \,\, \forall d_x \in {\cal C}(x^*).
\end{equation}
From the expression of  $\nabla^2 \varphi (x)$ in (\ref{phi-2deri}), we obtain (\ref{eq:secNCs}) from
(\ref{eq:SONC}). The proof is complete. \hfill $\Box$
\begin{theorem}\label{th:Sc}
(Second-order Sufficient Optimality Conditions)
Let $(x^*,y^*) \in \Re^n \times \Re^m$ be a point around which $f,h,g$ are twice continuously differentiable and $H$, $G$ are twice continuously differentiable around $x^*$. Assume that $x^* \in \Phi$ and $y^* \in Y(x^*)$.  Let $(\mu^*,\lambda^*) \in \Re^{m_1}\times \Re^{m_2}$.
 Suppose that  Problem $({\rm P}_{x^*})$  satisfies  Jacobian uniqueness conditions at $(y^*,\mu^*,\lambda^*)$, $\Lambda (x^*) \ne \emptyset$, and for every $d_x \in {\cal C}(x^*)\setminus \emptyset$ (where ${\cal C}(x^*)$ is defined by (\ref{eq:CritiCS})),
\begin{equation}\label{eq:secSCs}
\begin{array}{l}
\displaystyle \sup_{(u,v) \in \Lambda(x^*)} \left\{ \left\langle \left[\displaystyle \sum_{j=1}^{n_1}u_i\nabla^2_{xx}H_j(x^*)+\displaystyle \sum_{i=1}^{n_2} v_i\nabla^2_{xx}G_i(x^*)\right]d_x,d_x \right\rangle\right\}\\[16pt]
 \quad \quad \quad \, +\left\langle \left[ \nabla^2_{xx}{\cal L}(x^*;y^*,\mu^*,\lambda^*)- N(x^*)^TK(x^*)^{-1}N(x^*)\right]d_x,d_x \right\rangle > 0,
 \end{array}
\end{equation}
where $K(x)$ is defined by (\ref{eq:Kx}) and $N(x)$ is defined by (\ref{Nx}).
Then there exist $\delta_1 \in (0,\delta_0)$, $\varepsilon_1 \in (0,\varepsilon_0)$ (where $\delta_0$ and $\varepsilon_0$ are given by Lemma (\ref{lem:imp}))and $\gamma_1>0$,$\gamma_2>0$ such that
for $x \in \textbf{B}_{\delta_1}(x^*)\cap \Phi$ and $y \in \textbf{B}_{\varepsilon_1}(y^*)\cap Y(x^*)$,
\begin{equation}\label{eq:2ndG}
f(x^*,y)+\gamma_1 \|y-y^*\|^2/2 \leq f(x^*,y^*)\leq \displaystyle \sup_{z \in Y(x) \cap \textbf{B}_{\varepsilon_0}(y^*)} f(x,z)-\gamma_2\|x-x^*\|^2/2,
\end{equation}
which implies that $(x^*,y^*)$ a local minimax point of Problem (\ref{cminimax}).
\end{theorem}
{\bf Proof}. As Jacobian uniqueness conditions hold at $y^*$ for Problem (${\rm P}_{x^*}$), we know that the local second-order descent condition holds for Problem (${\rm P}_{x^*}$) at $y^*$. Thus there exist $\gamma_1>0$ and $\varepsilon_1\in (0,\varepsilon_0)$ such that
$$
f(x^*,y)+\gamma_1 \|y-y^*\|^2/2 \leq f(x^*,y^*),\ y \in \textbf{B}_{\varepsilon_1}(y^*)\cap Y(x^*).
$$
From the formula for $\nabla \varphi (x)$ in (\ref{phi-1deri}) and the formula for $\nabla^2 \varphi (x)$ in (\ref{phi-2deri}),
 we have from the definition (\ref{eq:CritiCS}) that
 $$
 {\cal C}(x^*)=\{d_x \in \Re^n: {\cal J}H(x^*)d_x=0; \nabla G_i(x^*)^Td_x \leq 0, i \in I(x^*); \nabla \varphi (x^*) ^Td_x\leq 0\},
 $$
  which is just the critical cone of Problem (\ref{prob:RP}), and (\ref{eq:secSCs})  is the second-order sufficient optimality condition for Problem (\ref{prob:RP}). Then the second-order growth condition of Problem (\ref{prob:RP}) holds at $x^*$ from \cite{BS00}; namely, there exist $\gamma_2 >$ and $\delta_1 \in (0, \delta_0)$  such that
  $$
  \varphi (x^*) + \gamma_2\|x-x^*\|^2/2\leq \varphi (x), \, x \in \textbf{B}_{\delta_1}(x^*)\cap \Phi,
  $$
 which combining  the expression
  $$
  \varphi (x)=\displaystyle \sup_{z \in Y(x) \cap \textbf{B}_{\varepsilon_0}(y^*)}f(x,z)
  $$
  yields
  $$
  f(x^*,y^*)\leq \displaystyle \sup_{z \in Y(x) \cap \textbf{B}_{\varepsilon_0}(y^*)} f(x,z)-\gamma_2\|x-x^*\|^2/2.
  $$
Therefore the inequalities in (\ref{eq:2ndG}) are demonstrated.
 \hfill $\Box$

 In the following, we
 derive necessary optimality conditions  for Problem (\ref{cminimax}) under {\bf Assumption A} for (${\rm P}_{x^*}$).
 Define the outer approximation of $\partial \varphi (x)$ by
\begin{equation}\label{eq:OS}
\begin{array}{l}
\widetilde\partial\varphi(x)=\nabla_x{\cal L}(x;y(x),\mu(x),\lambda(x))\\
\quad  \,-\left\{
 H(x,W)^T\nabla_{y,\mu,\lambda}{\cal L}(x;y(x),\mu(x),\lambda(x))
: W \in \partial \Pi_{\Re^{m_2}_-}(\lambda(x)+g(x,y(x)))
\right \}.
\end{array}
\end{equation}
\begin{theorem}\label{th:nc2}
(First-order Necessary Optimality Conditions \rev{under {\bf Assumption A}}) Let $(x^*,y^*) \in \Re^n \times \Re^m$ be a point around which $f,h,g$ are twice continuously differentiable and $H$, $G$ are  continuously differentiable around $x^*$. Let $(x^*,y^*)$  be a local minimax point of Problem (\ref{cminimax}). Assume  that linear independence constraint qualification holds at $y^*$ for constraint set $Y(x^*)$. Then there exists a unique vector $(\mu^*,\lambda^*) \in \Re^{m_1}\times\Re^{m_2}$ such that
\begin{equation}\label{KKT-Px2}
\begin{array}{l}
\nabla_y {\cal L}(x^*;y^*,\mu^*,\lambda^*)=0,\\[3pt]
h(x^*,y^*)=0,\\[3pt]
0\geq \lambda^* \perp g(x^*,y^*) \leq 0.
\end{array}
\end{equation}
For any $d_y \in {\cal C}_{x^*}(y^*)$, we have that
\begin{equation}\label{second-N2}
\langle \nabla^2_{yy}{\cal L}(x^*;y^*,\mu^*,\lambda^*)d_y, d_y \rangle \leq 0.
\end{equation}
 Suppose that  Problem $({\rm P}_{x^*})$  satisfies  {\bf Assumption A} at $(y^*,\mu^*,\lambda^*)$ and the
  Mangasarian-Fromovitz constraint qualification holds at $x^*$ for the constraint set $\Phi$. Then
 there exists $(u^*,v^*) \in \Re^{n_1}\times \Re^{n_2}$ and
 $W \in \partial \Pi_{\Re^{m_2}_-}(\lambda^*+g(x^*,y^*))$ such that
 \begin{equation}\label{KKT-P2}
\begin{array}{l}
\nabla_x {\cal L}(x^*;y^*,\mu^*,\lambda^*)-H(x^*,W)^T\nabla_{y,\mu,\lambda}{\cal L}(x^*;y^*,\mu^*,\lambda^*)+{\cal J}H(x^*)^Tu^*+{\cal J}G(x^*)^Tv^*=0,\\[3pt]
H(x^*)=0,\\[3pt]
0\leq v^* \perp G(x^*) \leq 0,
\end{array}
\end{equation}
where $H(x,W)$ is defined by (\ref{eq:Hxw}).
The set of all $(u^*,v^*)$ satisfying (\ref{KKT-P2}), denoted by $\Lambda (x^*)$, is a nonempty compact convex set.
\end{theorem}
{\bf Proof}. Properties (\ref{KKT-Px2}) and (\ref{second-N2}) are obvious from Theorem \ref{th:nc}. Now we prove property
(\ref{KKT-P2}).  It  follows from Corollary \ref{phi-d} that $\varphi$ is locally Lipschitz continuous, also directionally differentiable in $\textbf{B}(x^*,\delta_0)$. Thus we can easily get that $0 \in \Re^n$ is an optimal solution to the following problem
\begin{equation}\label{eq:ACLP}
\begin{array}{cl}
\min_{d_x} & \varphi'(x^*;d_x)\\[4pt]
{\rm s.t.} &  d \in T_{\Phi}(x^*),
\end{array}
\end{equation}
where $T_{\Phi}(x^*)$ is the tangent cone of $\Phi$ at $x^*$. Since the Mangasarian-Fromovitz constraint qualification holds at $x^*$ for the constraint set $\Phi$, we have from \cite{BS00} that
\begin{equation}\label{eq:TC}
T_{\Phi}(x^*)=\{d_x \in \Re^n: {\cal J}H(x^*)d_x =0, \nabla G_i(x^*)^Td_x \leq 0, i \in I(x^*)\}.
\end{equation}
Since $\partial \varphi (x^*) \subset \widetilde \partial \varphi (x^*)$, from (\ref{eq:objD}) and (\ref{eq:objCD}), we obtain
\begin{equation}\label{eq:OA}
\varphi'(x^*;d_x) \leq \max\{v^Td_x: v \in \widetilde \varphi (x^*)\}=\delta^*(d_x\,|\, \widetilde \varphi (x^*)).
\end{equation}
It follows from (\ref{eq:TC}) and (\ref{eq:OA}) and $0\in \Re^n$ is the minimizer of Problem (\ref{eq:ACLP}), we have that $0\in \Re^n$ is the optimal solution to the following convex problem
\begin{equation}\label{eq:ACLPT}
\begin{array}{cl}
\min_{d_x} & \delta^*(d_x\,|\, \widetilde \partial \varphi (x^*))\\[4pt]
{\rm s.t.} &  \nabla H_j(x^*)^Td_x=0,\, j=1,\ldots,n_1,\\[4pt]
& \nabla G_i(x^*)^Td_x=0,\, i \in I(x^*).
\end{array}
\end{equation}
Noting that the  Mangasarian-Fromovitz constraint qualification holds at $x^*$ for the constraint set $\Phi$, we have that Slater condition holds for convex optimization problem (\ref{eq:ACLPT}).
Then, from the optimality conditions for convex programming, we have that  there exist $u^*$ and $v^*_i, i \in I(x^*)$ such that
\begin{equation}\label{Conv-KKT}
0\in \widetilde \partial \varphi (x^*)+{\cal J}H(x^*)^Tu^*+\displaystyle \sum_{i \in I(x^*)} v^*_i\nabla G_i(x^*).
\end{equation}
Therefore
 there exists
 $W \in \partial \Pi_{\Re^{m_2}_-}(\lambda^*+g(x^*,y^*))$ such that
 (\ref{KKT-P2}) is satisfied. Now we prove by contradiction that $\Lambda (x^*)$ is compact. Assuming that the set is unbounded, there there exist a sequence $W^k \partial  \Pi_{\Re^{m_2}_-}(\lambda^*+g(x^*,y^*))$, $u^k$ and $v^k$ such that $(W^k,u^k,v^k)$ satisfies (\ref{KKT-P2}) and $\|(u^k,v^k)\| \rightarrow \infty$.  Let $(\bar u^k,\bar v^k)=(u^k,v^k)/\|(u^k,v^k)\|$ and without loss of generality assume that $(\bar u^k,\bar v^k)\rightarrow (\bar u,\bar v)$.  Thus we have
\begin{equation}\label{KKT-Pk}
\begin{array}{l}
\nabla_x {\cal L}(x^*;y^*,\mu^*,\lambda^*)/\|(u^k,v^k)\|+{\cal J}H(x^*)^T\bar u^k+{\cal J}G(x^*)^T\bar v^k=0,\\[3pt]
H(x^*)=0,\\[3pt]
0\leq \bar v^k \perp G(x^*) \leq 0.
\end{array}
\end{equation}
Taking the limit of (\ref{KKT-Pk}) as $k \rightarrow \infty$, we obtain
\begin{equation}\label{KKT-Pkbar}
\begin{array}{l}
{\cal J}H(x^*)^T\bar u+{\cal J}G(x^*)^T\bar v=0,\\[3pt]
H(x^*)=0,\\[3pt]
0\leq  \bar v \perp G(x^*) \leq 0.
\end{array}
\end{equation}
The system (\ref{KKT-Pkbar}) implies $u=0$ and $v=0$ from the Mangasarian-Fromovitz constraint qualification, this contradicts with
$\|(\bar u,\bar v)\|=1$. The proof is complete. \hfill $\Box$

 \section{Some Concluding Remarks}\label{Sec4}
In this paper, we have successfully extended the definition of local minimax point 
from unconstrained minimax optimization problems.
to constrained minimax optimization problems. By analyzing Jacobian uniqueness conditions
for the lower-level maximization problem and the strong regularity of Karush-Kuhn-Tucker conditions of the
maximization problem, we provided both necessary optimality conditions and sufficient optimality conditions
for the local minimax points of constrained minimax optimization problems, see Theorems \ref{th:nc}-\ref{th:nc2}.

As the current study is theoretical, we are looking for more applications of the constrained minimax optimization
problems in modern machine learning and also traditional research areas in future. Furthermore, it remains under
investigation how to design numerical algorithms which can converge to a local minimax point of constrained
minimax optimization problems.


\end{document}